\documentstyle[a4,twoside,amssymb]{article}

\catcode`\@=11
\long\def\@makefntext#1{
\protect\noindent \hbox to 3.2pt {\hskip-.9pt  
$^{{\eightrm\@thefnmark}}$\hfil}#1\hfill}		

\def\@makefnmark{\hbox to 0pt{$^{\@thefnmark}$\hss}}	
	
\def\ps@myheadings{\let\@mkboth\@gobbletwo
\def\@oddhead{\hbox{}
\rightmark\hfil\eightrm\thepage}   
\def\@oddfoot{}\def\@evenhead{\eightrm\thepage\hfil
\leftmark\hbox{}}\def\@evenfoot{}
\def\sectionmark##1{}\def\subsectionmark##1{}}

\oddsidemargin=\evensidemargin
\newcounter{sectionc}\newcounter{subsectionc}\newcounter{subsubsectionc}
\renewcommand{\section}[1] {\vspace{12pt}\addtocounter{sectionc}{1} 
\setcounter{subsectionc}{0}\setcounter{subsubsectionc}{0}\noindent 
	{\tenbf\thesectionc. #1}\par\vspace{5pt}}
\renewcommand{\subsection}[1] {\vspace{12pt}\addtocounter{subsectionc}{1} 
	\setcounter{subsubsectionc}{0}\noindent 
	{\bf\thesectionc.\thesubsectionc. {\kern1pt \bfit #1}}\par\vspace{5pt}}
\renewcommand{\subsubsection}[1] {\vspace{12pt}\addtocounter{subsubsectionc}{1}
	\noindent{\tenrm\thesectionc.\thesubsectionc.\thesubsubsectionc.
	{\kern1pt \tenit #1}}\par\vspace{5pt}}
\newcommand{\nonumsection}[1] {\vspace{12pt}\noindent{\tenbf #1}
	\par\vspace{5pt}}

\newcounter{appendixc}
\newcounter{subappendixc}[appendixc]
\newcounter{subsubappendixc}[subappendixc]
\renewcommand{\thesubappendixc}{\Alph{appendixc}.\arabic{subappendixc}}
\renewcommand{\thesubsubappendixc}
	{\Alph{appendixc}.\arabic{subappendixc}.\arabic{subsubappendixc}}

\renewcommand{\appendix}[1] {\vspace{12pt}
        \refstepcounter{appendixc}
        \setcounter{figure}{0}
        \setcounter{table}{0}
        \setcounter{lemma}{0}
        \setcounter{theorem}{0}
        \setcounter{corollary}{0}
        \setcounter{definition}{0}
        \setcounter{equation}{0}
        \renewcommand{\thefigure}{\Alph{appendixc}.\arabic{figure}}
        \renewcommand{\thetable}{\Alph{appendixc}.\arabic{table}}
        \renewcommand{\theappendixc}{\Alph{appendixc}}
        \renewcommand{\thelemma}{\Alph{appendixc}.\arabic{lemma}}
        \renewcommand{\thetheorem}{\Alph{appendixc}.\arabic{theorem}}
        \renewcommand{\thedefinition}{\Alph{appendixc}.\arabic{definition}}
        \renewcommand{\thecorollary}{\Alph{appendixc}.\arabic{corollary}}
        \renewcommand{\theequation}{\Alph{appendixc}.\arabic{equation}}
        \noindent{\tenbf Appendix#1}\par\vspace{5pt}}
\newcommand{\subappendix}[1] {\vspace{12pt}
        \refstepcounter{subappendixc}
        \noindent{\bf Appendix \thesubappendixc. {\kern1pt \bfit #1}}
	\par\vspace{5pt}}
\newcommand{\subsubappendix}[1] {\vspace{12pt}
        \refstepcounter{subsubappendixc}
        \noindent{\rm Appendix \thesubsubappendixc. {\kern1pt \tenit #1}}
	\par\vspace{5pt}}

\newcommand{\textlineskip}{\baselineskip=13pt}
\newcommand{\smalllineskip}{\baselineskip=10pt}

\def\eightcirc{
\begin{picture}(0,0)
\put(4.4,1.8){\circle{6.5}}
\end{picture}}
\def\eightcopyright{\eightcirc\kern2.7pt\hbox{\eightrm c}} 

\newcommand{\copyrightheading}[1]
	{\vspace*{-2.5cm}\smalllineskip{\flushleft
	{\footnotesize Mathematical Models and Methods in Applied Sciences #1}\\
	{\footnotesize $\eightcopyright$\, World Scientific Publishing
	 Company}\\
	 }}

\def\abstracts#1#2#3{{
	\centering{\begin{minipage}{4.5in}\baselineskip=10pt\footnotesize
	\parindent=0pt #1\par 
	\parindent=15pt #2\par
	\parindent=15pt #3
	\end{minipage}}\par}} 

\def\keywords#1{{
	\centering{\begin{minipage}{4.5in}\baselineskip=10pt\footnotesize
	{\footnotesize\it Keywords}\/: #1
	 \end{minipage}}\par}}

\renewenvironment{thebibliography}[1]
	{\frenchspacing
	 \ninerm\baselineskip=11pt
	 \begin{list}{\arabic{enumi}.}
        {\usecounter{enumi}\setlength{\parsep}{0pt}     
	 \setlength{\leftmargin 12.7pt}{\rightmargin 0pt} 
         \setlength{\itemsep}{0pt} \settowidth
	{\labelwidth}{#1.}\sloppy}}{\end{list}}

\newcounter{itemlistc}
\newcounter{romanlistc}
\newcounter{alphlistc}
\newcounter{arabiclistc}

\newcommand{\fcaption}[1]{
        \refstepcounter{figure}
        \setbox\@tempboxa = \hbox{\footnotesize Fig.~\thefigure. #1}
        \ifdim \wd\@tempboxa > 5in
           {\begin{center}
        \parbox{5in}{\footnotesize\smalllineskip Fig.~\thefigure. #1}
            \end{center}}
        \else
             {\begin{center}
             {\footnotesize Fig.~\thefigure. #1}
              \end{center}}
        \fi}

\newcommand{\tcaption}[1]{
        \refstepcounter{table}
        \setbox\@tempboxa = \hbox{\footnotesize Table~\thetable. #1}
        \ifdim \wd\@tempboxa > 5in
           {\begin{center}
        \parbox{5in}{\footnotesize\smalllineskip Table~\thetable. #1}
            \end{center}}
        \else
             {\begin{center}
             {\footnotesize Table~\thetable. #1}
              \end{center}}
        \fi}

\def\@citex[#1]#2{\if@filesw\immediate\write\@auxout
	{\string\citation{#2}}\fi
\def\@citea{}\@cite{\@for\@citeb:=#2\do
	{\@citea\def\@citea{,}\@ifundefined
	{b@\@citeb}{{\bf ?}\@warning
	{Citation `\@citeb' on page \thepage \space undefined}}
	{\csname b@\@citeb\endcsname}}}{#1}}

\newif\if@cghi
\def\cite{\@cghitrue\@ifnextchar [{\@tempswatrue
	\@citex}{\@tempswafalse\@citex[]}}
\def\citelow{\@cghifalse\@ifnextchar [{\@tempswatrue
	\@citex}{\@tempswafalse\@citex[]}}
\def\@cite#1#2{{$\null^{#1}$\if@tempswa\typeout
	{IJCGA warning: optional citation argument 
	ignored: `#2'} \fi}}

\def\pmb#1{\setbox0=\hbox{#1}
	\kern-.025em\copy0\kern-\wd0
	\kern.05em\copy0\kern-\wd0
	\kern-.025em\raise.0433em\box0}

\def\fnt#1#2{\footnotetext{\kern-.3em
	{$^{\mbox{\scriptsize #1}}$}{#2}}}

\def\fpage#1{\begingroup
\voffset=.3in
\thispagestyle{empty}\begin{table}[b]\centerline{\footnotesize #1}
	\end{table}\endgroup}

\def\runninghead#1#2{\pagestyle{myheadings}
\markboth{{\protect\footnotesize\it{\quad #1}}\hfill}
{\hfill{\protect\footnotesize\it{#2\quad}}}}
\font\tenrm=cmr10
\font\tenit=cmti10 
\font\tenbf=cmbx10
\font\bfit=cmbxti10 at 10pt
\font\ninerm=cmr9

\font\eightrm=cmr8

\newtheorem{theorem}{Theorem}

\newtheorem{lemma}{Lemma}

\def\qed{\hbox{${\vcenter{\vbox{			
   \hrule height 0.4pt\hbox{\vrule width 0.4pt height 6pt
   \kern5pt\vrule width 0.4pt}\hrule height 0.4pt}}}$}}

\def\theequation{\thesectionc.\arabic{equation}}	

\begin{document}

\newcommand{\pe}{\psi}
\def\d{\delta} 
\def\ds{\displaystyle} 
\def\e{{\epsilon}} 
\def\eb{\bar{\eta}}  
\def\enorm#1{\|#1\|_2} 
\def\Fp{F^\prime}  
\def\Kc{{\cal K}}
\def\norm#1{\|#1\|} 
\def\wb{{\bar w}} 
\def\zb{{\bar z}} 

\def\Re{{\Bbb R}}
\def\Ze{{\Bbb Z}}
\def\Na{{\Bbb N}}
\def\charf {{{\rm 1}\kern-.24em {\rm l}}} 
\newcommand{\Ind}[1]{\mbox{\Large \bf 1}_{#1}} 
\newcommand{\Ph}{{\cal P}^{h}}
\def\DT{\Delta t}
\def\undemi{\frac{1}{2}}
\def\sgn{\mbox{sgn}}
\def\ds{\displaystyle} 
\def\Dp#1#2{\partial_{#1}{#2}}
\def\ae{a^{\e}}
\def\be{b^{\e}}
\def\ue{u^{\e}}
\def\tue{\tilde{u}^{\e}}
\def\ve{v^{\e}}
\def\tve{\tilde{v}^{\e}}
\newtheorem{defi}{Definition}
\newtheorem{coro}{Corollary}
\newtheorem{proposition}{Proposition}
\newtheorem{remark}{Remark}
\def\BE{\begin{equation}}
\def\EE{\end{equation}}
\def \Dt{\partial_t}
\def \Dx{\partial_x}


\runninghead{Laurent Gosse}
{Measure source terms for a relaxing system}

\normalsize\textlineskip
\thispagestyle{empty}
\setcounter{page}{1}

\copyrightheading{}			

\vspace*{0.88truein}

\fpage{1}
\centerline{\bf TIME--SPLITTING SCHEMES AND MEASURE SOURCE TERMS}
\vspace*{0.035truein}
\centerline{\bf  FOR A QUASILINEAR RELAXING SYSTEM}
\vspace*{0.37truein}
\centerline{\footnotesize LAURENT GOSSE}
\vspace*{0.015truein}
\centerline{\footnotesize\it Istituto per le Applicazioni del Calcolo (sezione di Bari)}
\baselineskip=10pt
\centerline{\footnotesize\it via Amendola 122/I - 70126 Bari, ITALY.}
\centerline{\footnotesize\it Tel: +39 080 5530719. Fax: +39 080 5588235}
\centerline{\footnotesize\it E-mail: {\sf{l.gosse@area.ba.cnr.it}}}
\vspace*{10pt}

\vspace*{0.225truein}

\vspace*{0.21truein}
\abstracts{
Several singular limits are investigated in the context of a $2 \times 2$ 
system arising for instance in the modeling of chromatographic processes.
In particular, we focus on the case where the relaxation term and a $L^2$
projection operator are concentrated on a discrete lattice by means of
Dirac measures. This formulation allows to study more easily some
time-splitting numerical schemes.}{}{}

\vspace*{10pt}


\vspace*{1pt}\textlineskip	
\vspace*{-0.5pt}
\noindent












\section{Introduction}

Relaxation is a phenomenon which appears in a wide variety of physical 
situations. In gas dynamics, it occurs when the gas isn't in thermodynamic 
equilibrium. In elasticity, it is usually referred to as a fading memory 
mechanism. Hyperbolic conservation laws with relaxation are also served as 
discrete kinetic models. Relaxation approximations to scalar conservation laws can
be considered for purely numerical purposes too; general references are
{\it e.g.} \cite{gtz,fra,jpt1,jpt2,jx,lev,kat,kkm,lt,liu,nat,tv}.

In \cite{liu}, T.P. Liu studied the following $2 \times 2$ simple model
which captures the basic features of this physical process:
\BE\label{tpl}
\left\{\begin{array}{c}
\Dt u+\Dx f(u,v)=0,\\
\ds \Dt v +\Dx g(u,v)=\frac{1}{\varepsilon}(V(u)-v),
\end{array}\right. \qquad x \in \Re, t>0 \mbox{ with } \partial_v f(u,v)<0.
\EE
The numerical issue in such a system is to handle properly the relaxing process in
its infinite strength limit $\varepsilon \rightarrow 0$ on coarse
computational grids, see \cite{rpb,jin}. This requirement is twofold as
it includes some numerical
stability properties and a consistency with the expected
asymptotic behavior in the context of weak solutions. In order to derive
numerical schemes meeting these criteria, attention has been driven onto several
kinds of efficient treatments for the relaxation term: see \cite{areg,gt,gtz,jx,np}.
Indeed, robust and efficient approximations can be
derived based on ``time-splitting" algorithms where one alternates at every
time-step the resolution of the purely convective part and the handling of
the reaction process by any O.D.E. solver.

It can be sometimes shown rigorously that these numerical
schemes converge using {\it e.g.} the BV compactness
framework together with appropriate entropy consistency
properties. Our main goal here is to complete such a study
using slightly different arguments. Roughly speaking, we are about to consider both the
projection and the reaction term as being ``localized'' on a lattice; {\it formally},
this amounts to going from (\ref{tpl}) to
\BE\label{tpl2}
\left\{\begin{array}{c}
\ds \Dt u+\Dx f(u,v)=\sum_{n \in \Na_*}(\Ph(u)-u)\delta(t-n\DT),\\
\ds \Dt v +\Dx g(u,v)=\sum_{n \in \Na_*}\Big(\frac{\DT}{\varepsilon}(V(u)-v)
+(\Ph(v)-v)\Big)\delta(t-n\DT),
\end{array}\right.
\EE
where $\delta$ stands for the Dirac measure in $t=0$, $\DT>0$ is
a given time-step and $\Ph$ denotes the standard $L^2$ projection
operator on piecewise constant functions related to a mesh-size $h>0$.
Indeed, this new way to handle singular source terms in hyperbolic equations
can be found in \cite{bou,moc,gt,vas} in the context of kinetic equations
or semilinear relaxation to the one-dimensional conservation law. Such an approach
based on ``measure source terms" allows in particular to recover a fix
proposed in \cite{lev} to suppress spurious behavior in the stiff regime.
Indeed, the idea proposed in \cite{lev} for reactive chemical flows consists
in switching the order in which the reaction and projection steps appear
inside the Dirac masses, as we shall see in  \S3.2, \S3.3.

In this paper, we plan to  work out only the simpler case of
a quasilinear $2 \times 2$ system with relaxation previously considered
in \cite{rpb,tv,kur}, see (\ref{systeme}). We shall mainly rely on an
extensive use of some {\it nonconservative products}, \cite{hc,lft,ray},
in order to give a {\it rigorous} definition of the right-hand side, see
Propositions \ref{prop-proj} and \ref{prop-relax}. In return, our formulation entails direct
stability and convergence proofs exploiting mainly estimates easily proved
at the level of the original differential problem; this allows
to bypass some lengthy computations.

Hence, this work is organized as follows: in \S2, we establish some
compactness properties and study the nonconservative singular system obtained
as the relaxation and projection terms are concentrated inside Dirac masses, see
(\ref{systeme3}). We introduce a notion of {\sl entropy solution} well suited for this
kind of problems, see Definition \ref{gene-entro-soln}. Then,
in \S3, we present some applications to two types of time-splitting
numerical schemes for which  we prove convergence towards the ``equilibrium equation"
as the relaxation parameter diverges, see Theorems \ref{2} and \ref{4}.
 
\section{Study of a singular relaxing system}
\subsection{A quasilinear $2 \times 2$ model}

We are interested in the Cauchy problem for the following balance laws:
\BE
\label{systeme}
\left\{\begin{array}{c}
\Dt(u+v)+\Dx f(u)=0,\\
\Dt v =\mu \left(A(u)-v\right), \\
u(.,0)=u_0, v(.,0)=v_0.
\end{array}\right. \qquad x \in \Re, t>0
\EE
This system arises {\it e.g.} in chromatography, \cite{fra,nat,raa}. In 
this context, $u,v \in [0,1]$ stand for some species densities, one contained 
in a fluid flowing through a fixed bed and the other being absorbed 
by the material on the bed. This adsorption process is modeled by the 
right-hand side $\mu R(u,v)=\mu (A(u)-v)$, $\mu \gg 1$ where $1/\mu$ is the
relaxation time. This model has been extensively studied within the theory of
BV functions, \cite{daf,vol}, by the authors of \cite{stw,tv} (see
also \cite{luk,kur,nat,rpb}) under the following hypotheses:
\BE
\label{hypos}
\left\{\begin{array}{c}
f(0)=0, f'(u) \geq 0; \\
A(0)=0, A(1)=1, A'(u) \geq 0. \\
\end{array}\right.
\EE
We recall here some of their main results.

\begin{defi} \label{entro-soln}
We say that the pair $(u,v)$ is an {\bf entropy solution} of (\ref{systeme}) if:
\begin{itemize}
\item[(i)] $u(.,t)$ and $v(.,t)$ lie in $L^1 \cap BV(\Re)$ for any
$t \in \Re^+$ and satisfy (\ref{systeme}) in the sense of distributions;
\item[(ii)] there exists $M \in \Re^+$ such that for any $t \geq 0, s \geq 0$:
$$
\|u(.,t)-u(.,s)\|_{L^1(\Re)}+\|v(.,t)-v(.,s)\|_{L^1(\Re)} \leq M |t-s|;
$$
\item[(iii)] there holds for any $k,l \in \Re$ and any nonnegative
test-function in ${\cal D}(\Re \times \Re_*^+)$:
\BE \label{zaza}
\begin{array}{c}
\Dt \Big [|u-k|+|v-l|\Big] + \Dx \Big |f(u)-f(k) \Big | \\
\leq \mu R(u,v)\Big [\sgn(v-l)-\sgn(u-k)\Big ].
\end{array}
\EE
\end{itemize}
\end{defi}

The relevance of this notion comes from the forthcoming theorems whose detailed
proofs are to be found in \cite{tv}. We stress however that a more general
definition of entropy solution to (\ref{systeme}) has been proposed within
the class of $L^{\infty}$ functions, \cite{nat}.

\begin{theorem} \label{EU-thm}
Let $u_0,v_0 \in L^1 \cap BV(\Re)$; under the hypotheses (\ref{hypos}), there
exists a unique entropy solution to (\ref{systeme}) which satisfies furthermore
for any $t \in \Re^+$:
\begin{itemize}
\item[(i)] $\mu  \|R(u,v)(.,t)\|_{L^1(\Re)} \leq M$ for some $M \in \Re^+$;
\item[(ii)] Given another set of initial data $\tilde{u}_0,\tilde{v}_0$ in $L^1
\cap BV(\Re)$, there holds:
$$
\|u(.,t)-\tilde{u}(.,t)\|_{L^1(\Re)}+\|v(.,t)-\tilde{v}(.,t)\|_{L^1(\Re)} \leq
\|u_0-\tilde{u}_0\|_{L^1(\Re)}+\|v_0-\tilde{v}_0\|_{L^1(\Re)}.
$$
\end{itemize}
\end{theorem}

The theory of Kru\v{z}kov, \cite{kru}, ensures that there exists a unique
entropy solution $w$ to the following ``equilibrium" equation since $A' \geq 0$:
\BE
\label{equilib}
\left\{\begin{array}{l}
\Dt (w + A(w)) +\Dx f(w) =0, \qquad x \in \Re, t>0\\
w(.,0)=w_0  \in L^1 \cap BV(\Re).
\end{array}\right.
\EE

\begin{theorem} \label{THM-2.2}
Under the assumptions of Theorem \ref{EU-thm} and if $w_0=u_0$,
there exists a constant $M' \in \Re^+$ such that:
$$
\|u(.,t)-w(.,t)\|_{L^1(\Re)} \leq \frac{M'}{\mu^{\frac{1}{3}}}, \quad t \in \Re^+.
$$
\end{theorem}

\subsection{Introduction of a ``localized" singular system}

In all the sequel, $\DT$ and $h$ will stand for fixed positive parameters
related to some cartesian computational discretization.

\begin{itemize}
\item We first introduce a projector onto piecewise constant functions; for any
$j \in \Ze$, we denote $C_j \stackrel{def}{=} [(j-\undemi)h, (j+\undemi)h[$.
Then, we define:
\BE \label{projector}
\begin{array}{cccl}
\Ph: & L^1 \cap BV(\Re) & \rightarrow & L^1 \cap BV(\Re) \\
&\varphi & \mapsto & \ds \left(\frac{1}{h}\int_{C_j}\varphi(x).dx \right)_{j \in \Ze}
\end{array}
\EE
\item Let $I^n \stackrel{def}{=} [n\DT, (n+1)\DT[$; we define two functions
which are Lipschitz continuous with respect to time for $\e>0$ small enough:
\BE \label{funct-a}
\Re^+ \ni t \mapsto \ae(t)=\left\{
\begin{array}{cl}
(n-1)\DT, & t \in [(n-1)\DT, n\DT- \e [ \\
\ds n\DT + \left(\frac{t -n\DT}{\e}\right)\DT, & t \in [n\DT- \e, n\DT[
\end{array}\right.
\EE
and:
\BE \label{funct-b}
\Re^+ \ni t \mapsto \be(t)=\left\{
\begin{array}{cl}
\ds n\DT + \left(\frac{t -n\DT}{\e}\right)\DT, & t \in [n\DT , n\DT+ \e[ \\
(n+1)\DT, & t \in [n\DT+ \e, (n+1)\DT [
\end{array}\right.
\EE
\end{itemize}

We can therefore consider the following system for any $\e >0$, $\mu,\nu$ in
$\Re^+$ within the theory proposed in \cite{tv}:
\BE
\label{systeme2}
\left\{\begin{array}{l}
\Dt \ue +\Dx f(\ue)=\nu (\Ph(\ue)-\ue)\Dt \ae - \mu R(\ue,\ve)\Dt \be, \\
\Dt \ve =\nu (\Ph(\ve)-\ve)\Dt \ae + \mu R(\ue,\ve)\Dt \be, \\
\ue(.,0)=u_0, \ve(.,0)=v_0.
\end{array}\right.
x \in \Re, t>0
\EE

In the limit $\e \rightarrow 0$, both $\ae$ and $\be$ converge almost
everywhere towards the piecewise constant function $t \mapsto \DT
\sum_{n \in \Na_*} Y(t -n\DT)$, where $Y$ stands for the classical Heaviside
distribution.

We give at once a relaxation estimate for (\ref{systeme2}) which
is uniform in $\e$.

\begin{lemma} \label{relaxing}
Under the assumptions of Theorem \ref{EU-thm}, the entropy solution of
(\ref{systeme2}) satisfies the relaxation estimate
$\mu \|R(\ue,\ve)\Dt{\be}\|_{{\cal M}_{loc}(\Re \times \Re_*^+)} \leq \bar{C}$
for some $\bar{C} \in \Re^+$.
\end{lemma}

{\bf Proof. }
We follow the classical ideas from \cite{kat,nat,tv} introducing a special
entropy/entropy-flux pair for (\ref{systeme2}):
$$
\eta(\ue,\ve)=\frac{(\ue)^2}{2}+H(\ve), \qquad H(\ve)=\int_0^{\ve}
A^{-1}(\xi).d\xi, \qquad q(\ue)=\int^{\ue}\xi.f'(\xi).d\xi.
$$
We observe that under the hypotheses (\ref{hypos}), the function $H$ is
strictly convex. Therefore, for any positive test-function, the following
inequality holds:
\BE \label{trou}
\begin{array}{c}
\Dt \eta(\ue,\ve)+\Dx q(\ue) \leq
-\mu \big(A(\ue)-\ve\big)\big(\ue-A^{-1}(\ve)\big)\Dt \be \\
+\nu \Big((\Ph(\ue)-\ue)\ue+(\Ph(\ve)-\ve)A^{-1}(\ve)\Big)\Dt \ae.
\end{array}
\EE
We integrate (\ref{trou}) on $\Re \times [0,T]$ for any $T \in \Re^+$.
Using Jensen's inequality, we can get rid of the terms involving the
projection operator since $\Dt \ae \geq 0$:
$$
\int_{\Re \times [0,T]} (\Ph(\ue)-\ue)\ue \Dt \ae.dx.dt
=h \int_{[0,T]} \Big(\Ph(\ue)^2-\Ph\big((\ue)^2\big)\Big) \Dt \ae.dt \leq 0.
$$
We observe now that since $A \in C^1(0,1)$ is increasing according
to (\ref{hypos}), the derivative of
its inverse mapping satisfies the following bound:
$$
\|(A^{-1})'\|_{C^0} \geq \frac{1}{\|A'\|_{C^0}} \stackrel{def}{=} C.
$$
And we have for any $w,v$ in $[0,1]$:
$$
\big(A^{-1}(w)-A^{-1}(v)\big)(w-v) \geq C (w-v)^2.
$$
Hence the integral involving $\Ph(\ve)$ is treated the following way:
$$
\begin{array}{c}
\ds \int_{\Re \times [0,T]} \big(\Ph(\ve)-\ve\big)A^{-1}(\ve) \Dt \ae.dx.dt= \\
-\ds \int_{\Re \times [0,T]} \big(\Ph(\ve)-\ve\big)
\big(A^{-1}(\Ph(\ve))-A^{-1}(\ve)\big) \Dt \ae.dx.dt \leq \\
\ds -C \int_{\Re \times [0,T]} \big(\Ph(\ve)-\ve\big)^2 \Dt \ae.dx.dt \leq 0.
\end{array}
$$
Finally, choosing $w=A(u)$, we have for $C>0$, the lower bound of
$(A^{-1})'$:
$$
(A(\ue)-\ve)(\ue-A^{-1}(\ve)) \geq C (A(\ue)-\ve)^2.
$$
Thus we get the expected relaxation estimate (see also \cite{kat}):
$$
C\mu \int_{\Re \times [0,T]} \big(A(\ue)-\ve\big)^2 \Dt \be .dx.dt \leq
\int_\Re \Big(\eta(\ue,\ve)(x,0)-\eta(\ue,\ve)\Big)(x,T).dx.
$$
Since $BV(\Re) \subset L^\infty(\Re)$, we are done.
$\Box$

Another important feature of (\ref{systeme2}) lies in the following compactness
result which will be of constant use in the sequel of the paper.

\begin{lemma}\label{compact}
Under the assumptions of Theorem \ref{EU-thm}, let $(\ue,\ve)$ be a sequence
of entropy solutions to (\ref{systeme2}). Then
$(\ue,\ve)$ is relatively
compact in $L^1_{loc}(\Re \times \Re^+_*)$ as $\e \rightarrow 0$.
\end{lemma}

{\bf Proof. }
It follows the classical BV stability canvas and makes use of the so-called
``quasi-monotonicity" property of the relaxation term, \cite{areg,nat,stw}.
\begin{itemize}

\item L$^1(\Re)$ bound: we multiply (\ref{systeme2}) by $(\sgn(\ue),\sgn(\ve))^T$
and integrate on $x \in \Re$. The point to be checked is the behavior of the
source terms, but we see that:
$$
\int_\Re \sgn(\ue)(\Ph(\ue)-\ue)(x,t).dx \leq
\sum_{j \in \Ze} \int_{C_j} |\ue(x,t)|.dx - \|\ue(.,t)\|_{L^1(\Re)} \leq 0.
$$
Also, we notice that $R(0,0)=0$ and we can use the mean-value theorem. Thanks
to the sign assumption on $A'$ in (\ref{hypos}), we get:
$$
\int_\Re(\sgn(\ve)-\sgn(\ue))(A'(\xi)\ue-\ve).dx \leq 0.
$$
Since $\Dt \ae$ and $\Dt \be$ are nonnegative, we derive finally:
\BE \label{Lun}
\|\ue(.,t)\|_{L^1(\Re)}+\|\ve(.,t)\|_{L^1(\Re)} \leq
\|u_0\|_{L^1(\Re)}+\|v_0\|_{L^1(\Re)}.
\EE

\item BV$(\Re)$ bound: Relying on \cite{kru}, we can differentiate the system
(\ref{systeme2}) with respect to
$x$, we multiply by $(\sgn(\Dx \ue),\sgn(\Dx \ve))^T$ and integrate on
$x \in \Re$. Using exactly the same arguments, we derive:
\BE \label{BeVe}
TV(\ue)(.,t)+TV(\ve)(.,t) \leq TV(u_0)+TV(v_0).
\EE

\item Let $t_0 \geq 0$ and $t=t_0 +\xi$, $\xi>0$; we compute
$$
\left\{\begin{array}{rcl}
\ds \int_\Re |\ue(x,t)-\ue(x,t_0)|.dx & \leq & 
\int_\Re \int_{t_0}^t |\Dx f(\ue)|.dt.dx  \\
&& +\int_\Re \int_{t_0}^t \mu|R(\ue,\ve)|\Dt \be .dt.dx\\
&& + \int_\Re \int_{t_0}^t \nu |\Ph(\ue)-\ue|\Dt \ae.dt.dx, \\

\ds \int_\Re |\ve(x,t)-\ve(x,t_0)|.dx & \leq & 
 \int_\Re \int_{t_0}^t \mu|R(\ue,\ve)|\Dt \be .dt.dx\\
&&+\int_\Re \int_{t_0}^t \nu |\Ph(\ve)-\ve|\Dt \ae.dt.dx, \\
\end{array}\right.
$$
in conjunction with the bounds:
$$
\left\{\begin{array}{l}
\ds \mu \int_{\Re \times [0,T]} |R(\ue,\ve)|\Dt \be .dx.dt \leq O(1), \\
\ds \int_\Re \Big(|\Ph(\ue)-\ue|+|\Ph(\ve)-\ve|\Big)(x,t).dx \leq h\big(
TV(u_0)+TV(v_0)\big).
\end{array}\right.
$$
And we derive:
\BE \label{equic}
\begin{array}{c}
\ds \sup_{\xi \not = 0} \int_{\Re \times [0,T]} 
\frac{|\ue(x,t+\xi)-\ue(x,t)|}{\xi}+\frac{|\ve(x,t+\xi)-\ve(x,t)|}{\xi}.dx.dt
\leq \\
\ds (T+\DT).\Big\{(\nu h+Lip(f))(TV(u_0)+TV(v_0)) +O(1)\Big\}
\end{array}
\EE
\end{itemize}
Therefore, we see that the sequence $(\ue,\ve)$ lies in BV$_{loc}(\Re \times \Re^+_*)$
and it remains to invoke Helly's compactness principle to conclude the proof. $\Box$

We notice that because of the Dirac masses in time arising at the right-hand side of
the equations, layers are likely to appear every time $n\DT$, $n \in \Na_*$ and
$$\ue,\ve \not \in C^0(\Re^+;L^1(\Re)),
\qquad \epsilon \rightarrow 0,$$
as in (for instance) \cite{pt} in the case of initial data being not ``well-prepared"; or it
can be evidenced just by taking $u_0 \equiv 0$, $v_0(x)=\sin(\pi x/h)$ inside (\ref{systeme3}).

\subsection{A meaning for the ambiguous products}

The preceding bounds for $\ue$ imply some compactness for the terms lying in 
the right-hand side in the weak-$\star$ topology of measures on $\Re \times 
\Re_*^+$. More precisely, we plan to shed some light on the products emanating 
in (\ref{systeme2}) in the limit $\e \rightarrow 0$. First, we state a stabilization result.

\begin{lemma} \label{lem-stab}
Under the assumptions of Theorem \ref{EU-thm}, the sequences of entropy solutions
$(u^\zeta,v^\zeta)$ to the Cauchy problem with initial data
$(u^\zeta_0,v^\zeta_0)\in L^1 \cap BV(\Re)$ for:
$$
\left.\begin{array}{c}
\Dt u^\zeta + \zeta.\Dx f(u^\zeta)=\nu \DT (\Ph(u^\zeta)-u^\zeta), \\
\Dt v^\zeta =\nu \DT (\Ph(v^\zeta)-v^\zeta), \\
\end{array}\right.
\mbox{or} \left.\begin{array}{c}
\Dt u^\zeta + \zeta.\Dx f(u^\zeta)=-\mu \DT R(u^\zeta,v^\zeta), \\
\Dt v^\zeta =\mu \DT R(u^\zeta,v^\zeta), \\
\end{array}\right.
$$
are relatively
compact in $L^1_{loc}(\Re \times \Re_*^+)$ and belong to $C^0(\Re^+;L^1(\Re))$
as $\zeta \rightarrow 0$.
\end{lemma}

{\bf Proof. }
It follows from the same arguments than the proof of Lemma \ref{compact}.
$\Box$

\begin{proposition} \label{prop-proj}
Under the hypotheses of Lemma \ref{compact}, $(\ue,\ve) \rightarrow (u,v)$
in $L^1_{loc}(\Re \times \Re_*^+)$ as $\e \rightarrow 0$ up to a subsequence.
Moreover, there holds:
\BE \label{lim-proj}
\begin{array}{c}
\displaystyle
\left(\begin{array}{c} \Ph(\ue)-\ue \\ \Ph(\ve)-\ve \end{array}\right)\Dt \ae 
\stackrel{weak-\star\ {\cal M}}{\rightharpoonup} \\
\displaystyle
\sum_{n \in \Na_*} \DT \left(\int_0^1 \left(\begin{array}{c} 
\Ph(\bar{u})-\bar{u} \\ \Ph(\bar{v})-\bar{v}\end{array}\right)(x,\tau)
d\tau\right) \delta(t-n\DT)
\end{array}
\EE
where $\bar{u},\bar{v}$ satisfy the following differential equations for 
$\tau \in [0,1]$:
\BE \label{edo-proj}
\partial_\tau \left(\begin{array}{c} \bar{u} \\ \bar{v} \end{array}\right)=
\nu \DT \left(\begin{array}{c} \Ph(\bar{u})-\bar{u} \\ \Ph(\bar{v})-\bar{v}
\end{array}\right),
\EE
together with the initial data, $x \in \Re, t=n\DT$ for $n \in \Na_*$:
\BE \label{init-proj}
\bar{u}(x,\tau=0)=u(x,t-0), \qquad \bar{v}(x,\tau=0)=u(x,t-0).
\EE
\end{proposition}

{\bf Proof. }
We want to compute the value of
the following expression for every $\phi \in C^0_c(\Re \times \Re_*^+)$, the
space of continuous functions with compact support in $\Re \times \Re^+_*$:
\BE \label{gros}
\begin{array}{c}
\ds \int_{\Re \times \Re_*^+} (\Ph(\ue)-\ue)\Dt \ae \phi(x,t).dx.dt = \\
\ds \sum_{n \in \Na_*} \int_\Re \int_{n\DT -\e}^{n\DT} (\Ph(\ue)-\ue)
\frac{\DT}{\e} \phi(x,t).dx.dt.
\end{array}
\EE
We pick a $n \in \Na_*$ and define $\tau=1+\frac{t-n\DT}{\e}$. Thus, 
according to this ``inner variable", the system (\ref{systeme2}) rewrites for 
$\tau \in [0,1]$:
$$
\left\{\begin{array}{c}
\partial_\tau \bar \ue +\e \Dx f(\bar \ue)=\nu \DT (\Ph(\bar \ue)- \bar \ue), \\
\partial_\tau \bar \ve =\nu \DT (\Ph(\bar \ve)-\bar \ve), \\
\end{array}\right.
$$
together with the initial data:
$$
\bar \ue(x,\tau=0)=\ue(x,n\DT-\e), \qquad \bar \ve(x,\tau=0)=\ve(x,n\DT-\e).
$$
We perform the change of variables $t \mapsto \tau$ in (\ref{gros}) to get:
$$
\begin{array}{c}
\ds \sum_{n \in \Na_*} \int_\Re \int_{n\DT -\e}^{n\DT} (\Ph(\ue)-\ue)
\frac{\DT}{\e} \phi(x,t).dx.dt=\\
\ds \sum_{n \in \Na_*} \int_\Re \int_{0}^{1} (\Ph(\bar \ue)-\bar \ue)(x,\tau)
\frac{\DT}{\e} \phi(x,n\DT+\e \tau -\e).dx.\e d\tau.
\end{array}
$$
Lemma \ref{lem-stab} implies that $\bar \ue$ satisfies (\ref{lim-proj}),
(\ref{edo-proj}) in the limit $\e \rightarrow 0$ and we are done.
$\Box$

\begin{remark}
Thanks to the linear form of the projection term, it is possible to define
directly the limit of $(\Ph(\ue)-\ue)\Dt \ae$ as $\e \rightarrow 0$ as a 
distribution of order zero ({\it i.e.} a bounded measure) provided $\DT=O(h)$. 
Let $\phi \in {\cal D}(\Re \times \Re^+_*)$, we can compute:
$$
\sum_{n \in \Na_*} \int_\Re (\Ph(u)-u)\phi(x,n\DT).dx=
\ds \sum_{j,n \in \Ze \times \Na_*} 
h \Big(\Ph(u)\Ph(\phi)-\Ph(u\phi)\Big)(jh,n\DT).
$$
And we have 
$\sum_{j \in \Ze}|\Ph(u)\Ph(\phi)-\Ph(u\phi)|(jh,t) \leq 
\|\phi\|_{C^0} TV(u)(.,t)$.
\end{remark}

Concerning the relaxation term, the situation could be different because of its 
nonlinear structure. Nevertheless, we have the following result.

\begin{proposition} \label{prop-relax}
Under the hypotheses of Lemma \ref{compact}, $(\ue,\ve) \rightarrow (u,v)$
in $L^1_{loc}(\Re \times \Re_*^+)$ as $\e \rightarrow 0$ up to a subsequence.
Moreover, there holds:
\BE \label{lim-relax}
\Big(A(\ue)-\ve\Big)\Dt \be \stackrel{weak-\star\ {\cal M}}{\rightharpoonup}
\sum_{n \in \Na_*} \DT \left(\int_0^1 \big(A(\bar{u})-\bar{v}\big)(x,\tau)
d\tau\right) \delta(t-n\DT)
\EE
where $\bar{u}, \bar{v}$ satisfy the following differential equations for 
$\tau \in [0,1]$:
\BE \label{edo-relax}
\partial_\tau \left(\begin{array}{c} \bar{u} \\ \bar{v} \end{array}\right)
=\mu \DT \left(\begin{array}{c} -\big(A(\bar{u})-\bar{v}\big) \\
\big(A(\bar{u})-\bar{v}\big) \end{array}\right),
\EE
together with the initial data, $x \in \Re, t=n\DT$ for $n \in \Na_*$:
\BE \label{init-relax}
\bar{u}(x,\tau=0)=u(x,t-0),\qquad \bar{v}(x,\tau=0)=v(x,t-0).
\EE
\end{proposition}

{\bf Proof. }
Let $\phi \in C^0_c(\Re \times \Re_*^+)$, we intend now to pass to the limit 
$\e \rightarrow 0$ in the following expression:
\BE \label{gros2}
\begin{array}{c}
\ds \int_{\Re \times \Re_*^+} \big(A(\ue)-\ve\big)\Dt \be \phi(x,t).dx.dt = \\
\ds \sum_{n \in \Na_*} \int_\Re \int_{n\DT}^{n\DT +\e} \big(A(\ue)-\ve\big)
\frac{\DT}{\e} \phi(x,t).dx.dt.
\end{array}
\EE
We pick a $n \in \Na_*$ and define $\tau=\frac{t-n\DT}{\e}$ as an ``inner 
variable" inside the system (\ref{systeme2}) which reads for $\tau \in [0,1]$:
$$
\left\{\begin{array}{c}
\partial_\tau \bar \ue +\e \Dx f(\bar \ue)=-\mu \DT \big(A(\bar \ue)-\bar \ve\big), \\
\partial_\tau \bar \ve =\mu \DT \big(A(\bar \ue)-\bar \ve\big), \\
\end{array}\right.
$$
together with the initial data:
$$
\bar \ue(x,\tau=0)=\ue(x,n\DT), \qquad \bar\ve(x,\tau=0)=\ve(x,n\DT).
$$
In order to conclude the proof, we perform the change of variables 
$t \mapsto \tau$ in (\ref{gros2}) and invoke Lemma \ref{lem-stab} 
to derive (\ref{lim-relax}), (\ref{edo-relax}) in the limit $\epsilon \rightarrow 0$.
$\Box$

\subsection{A contraction result ``\`a la Kru\v{z}kov"}

We can now consider weak solutions of the following singular relaxing system
for $\mu,\nu$ in $\Re^+$, $x \in \Re, t>0$:
\BE
\label{systeme3}
\left\{\begin{array}{c}
\Dt u +\Dx f(u)=\sum_{n \in \Na_*} \DT
\Big[\nu (\Ph(u)-u) - \mu R(u,v)\Big]\delta(t-n\DT), \\
\Dt v =\sum_{n \in \Na_*} \DT
\Big[\nu (\Ph(v)-v) + \mu R(u,v)\Big]\delta(t-n\DT), \\
u(.,0)=u_0 \in L^1 \cap BV(\Re), v(.,0)=v_0 \in L^1 \cap BV(\Re),
\end{array}\right.
\EE
where all the measure-valued source terms are {\bf always} to be understood in the sense
of distributions according  to the
results stated in Propositions \ref{prop-proj} and \ref{prop-relax}. We first
extend the notion of ``entropy solution" to (\ref{systeme3}).

\begin{defi} \label{gene-entro-soln}
We say that the pair $(u,v)$ is an {\bf entropy solution} of
(\ref{systeme3}) if:
\begin{itemize}
\item[(i)] $u$ and $v$ lie in $L^1 \cap BV(\Re \times [0,T])$ for any 
$T \in \Re^+$ and satisfy (\ref{systeme3}) in the sense of distributions;
\item[(ii)] there holds for any $k,l \in \Re$ and any nonnegative 
test-function in ${\cal D}(\Re \times \Re_*^+)$:
\BE \label{zozo}
\begin{array}{c}
\ds \Dt \Big [|u-k|+|v-l|\Big] + \Dx \Big |f(u)-f(k) \Big | \leq \\
\ds \sum_{n \in \Na_*} \DT \Big[\nu (\Ph(u)-u)\sgn(u-k)+\nu (\Ph(v)-v)\sgn(v-l)+ \\
\ds \mu R(u,v)\big (\sgn(v-l)-\sgn(u-k)\big )\Big]\delta(t-n\DT).\\
\end{array}
\EE
\end{itemize}
\end{defi}
The main result of this section states that system (\ref{systeme3}),
(\ref{zozo}) is $L^1(\Re)$-contractive.

\begin{theorem} \label{EU-sing}
Let $(u_0,v_0) \in L^1 \cap BV(\Re)$; under the hypotheses (\ref{hypos}), there
exists a unique entropy solution $(u,v)$ to (\ref{systeme3}). 
Moreover, if $(\tilde{u}_0,\tilde{v}_0)$ stands for another set of initial 
data in $L^1 \cap BV(\Re)$ and $(\tilde{u},\tilde{v})$ for its  associated 
entropy solution, there holds for any $t \in \Re^+$:
$$
\|u(.,t)-\tilde{u}(.,t)\|_{L^1(\Re)}+\|v(.,t)-\tilde{v}(.,t)\|_{L^1(\Re)} \leq 
\|u_0-\tilde{u}_0\|_{L^1(\Re)}+\|v_0-\tilde{v}_0\|_{L^1(\Re)}.
$$
\end{theorem}

Despite the fact that the solution of (\ref{systeme3}), (\ref{zozo}) suffer 
discontinuities on the lines $t=n\DT$, $n \in \Na_*$, its BV regularity allows 
one to define a set ${\cal N} \subset [0,T]$ of measure zero containing 
$0$ and $T$ such that each point in the strip $(x,t) \in \Re \times [0,T]$,
$t \not \in  {\cal N}$ is either an approximate jump or an approximate 
continuity location, \cite{daf}. 

{\bf Proof. }
We need first to check the regularity criterion (i) in Definition 
\ref{gene-entro-soln}. This is indeed a direct consequence of the uniform 
bounds (\ref{Lun}), (\ref{BeVe}), (\ref{equic}) in the proof of Lemma 
\ref{compact}. 
We start from the entropy inequality (\ref{zaza}) satisfied by $\ue,\ve$ for 
$\e>0$:
\BE \label{sta}
\begin{array}{c}
\ds \Dt \Big [|\ue-k|+|\ve-l|\Big] + \Dx \Big |f(\ue)-f(k) \Big | \leq \\
\ds \nu \big[(\Ph(\ue)-\ue)\sgn(\ue-k)+ (\Ph(\ve)-\ve)\sgn(\ve-l) \big] \Dt \ae+ \\
\ds \mu R(\ue,\ve)\big (\sgn(\ve-l)-\sgn(\ue-k)\big )\Dt \be.\\
\end{array}
\EE
Therefore, the inequality (\ref{zozo}) comes in the limit of (\ref{sta}) by
passing to the limit $\e \rightarrow 0$ as in the proofs of Propositions 
\ref{prop-proj} and \ref{prop-relax}. 

Concerning the contraction property which implies uniqueness, we follow the 
ideas of \cite{kru,kuz,tv}. Let us now introduce a nonnegative test function 
$\Psi(x,t,y,s)=\psi(x,t).\zeta(x-y).\zeta(t-s)$ belonging to ${\cal D}((\Re 
\times \Re_*^+)^2)$ with $\zeta$ the standard approximation of the Dirac mass.
We test the preceding entropy inequalities (\ref{sta}) for both $\ue,\ve$ and $\tue,\tve$
with this particular function, we add and integrate on $(\Re \times \Re_*^+)^2$.
For $\e>0$, the standard theory of Kru\v{z}kov allows to let $\zeta$ concentrate 
to the Dirac measure since $\ae, \be$ are Lipschitz functions. We obtain:
\BE \label{oups}
\begin{array}{c}
\ds \Dt \Big [|\ue-\tue|+|\ve-\tve|\Big] + \Dx \Big |f(\ue)-f(\tue) \Big | \leq \\
\ds \nu \big[((\Ph(\ue)-\ue)-(\Ph(\tue)-\tue))\sgn(\ue-\tue)+ \\
((\Ph(\ve)-\ve)-(\Ph(\tve)-\tve))\sgn(\ve-\tve) \big] \Dt \ae+ \\
\ds \mu (R(\ue,\ve)-R(\tue,\tve))\big (\sgn(\ve-\tve)-\sgn(\ue-\tue)\big )\Dt \be.\\
\end{array}
\EE

We study mainly the zero-order terms. Using the mean-value theorem for the 
relaxation term leads to
$$
\begin{array}{c}
\ds \int_{\Re \times \Re_*^+}\Big(\sgn(\ve-\tve)-\sgn(\ue-\tue)\Big)
\Big(A(\ue)-A(\tue)-\ve+\tve\Big)\Dt \be \psi(x,t).dx.dt \\
\ds =\int_{\Re \times \Re_*^+}A'(\xi) 
\Big(\sgn(\ve-\tve)(\ue-\tue)-|\ue-\tue|\Big)\Dt \be \psi(x,t).dx.dt \\
+\ds \int_{\Re \times \Re_*^+}
\Big(\sgn(\ue-\tue)(\ve-\tve)-|\ve-\tve|\Big)\Dt \be \psi(x,t).dx.dt,
\end{array}
$$
for a $\xi \in [\min(\ue,\tue),\max(\ue,\tue)]$. Since $\Dt \be \geq 0$ and 
$A' \geq 0$, this term gives a nonpositive contribution on the right-hand side 
of the entropy inequality (\ref{oups}). Concerning the projection operator, we 
see that
$$
\begin{array}{c}
\ds \int_{\Re\times \Re_*^+} \sgn(\ue-\tue)
\Big(\Ph(\ue-\tue)-(\ue-\tue)\Big)\Dt \ae \psi(x,t).dx.dt = \\
\ds \int_{\Re \times \Re_*^+}\Big(\sgn(\ue-\tue)-\sgn(\Ph(\ue-\tue))\Big)
\Big(\Ph(\ue-\tue)-(\ue-\tue)\Big)\Dt \ae \psi  .dx.dt\\
\ds +\int_{\Re \times \Re_*^+}\sgn\big(\Ph(\ue-\tue)\big)
\Big(\Ph(\ue-\tue)-(\ue-\tue)\Big)\Dt \ae \psi(x,t) .dx.dt.
\end{array}
$$
Since $\Ph(\ue-\tue)$ is constant on each cell $C_j$, the last integral 
vanishes. Once again, since $\Dt \ae \geq0$, this term gives a nonpositive
contribution in (\ref{oups}). We end up with:
$$
\Dt \Big [|\ue-\tue|+|\ve-\tve|\Big] + \Dx \Big |f(\ue)-f(\tue) \Big | \leq 0.
$$

It remains now to select $\psi$ as in \cite{kru,kuz,bp} to complete the proof 
of Theorem \ref{EU-sing}. $\Box$

\begin{remark}
We can notice that this contraction property gives back the $BV(\Re)$
regularity for $u(.,t)$, $v(.,t)$ as soon as $u_0,v_0$ belong to $BV(\Re)$ 
since the problem (\ref{systeme3}) is translation invariant.
\end{remark}

We notice finally that since $\Ph$ is linear, any weak solution of 
(\ref{systeme3}) also satisfies:
$$
\left\{\begin{array}{c}
\Dt (u+v) +\Dx f(u)=\sum_{n \in \Na_*} \nu \DT
\Big[\Ph(u+v)-(u+v)\Big]\delta(t-n\DT), \\
\Dt v =\sum_{n \in \Na_*} \DT
\Big[\nu (\Ph(v)-v) + \mu R(u,v)\Big]\delta(t-n\DT). \\
\end{array}\right.
$$

\section{Study of the resulting ``time--splitting" numerical schemes}
\subsection{A technical lemma}

We point out quickly that the special structure of the ``singular projection 
term" allows us to refine a bit the entropy inequality (\ref{zozo}). 

\begin{lemma} \label{entro-proj}
In the sense of Proposition \ref{prop-proj}, the following holds for 
any $k \in \Re$ and any nonnegative continuous compactly supported 
test-function:
$$
\begin{array}{c}
\ds \sum_{n \in \Na_*} \nu\DT \Big(\Ph(u)-u\Big)\sgn(u-k).\delta(t-n\DT) \leq\\
\ds \sum_{n \in \Na_*}\Big(|\Ph(u)-k|-|u-k|\Big) \big(1-\exp(-\nu \DT)\big).
\delta(t-n\DT).
\end{array}
$$
\end{lemma}

{\bf Proof. }
Let $0 \leq \phi \in C^0_c(\Re \times \Re_*^+)$; we perform the following 
splitting:
$$
\left\{\begin{array}{l}
\sgn(\bar{u}-k)=\sgn(\Ph(\bar{u})-k)+\big(\sgn(\bar{u}-k)-\sgn(\Ph(\bar{u})-k)\big), \\
\Ph(\bar{u})-\bar{u}=(\Ph(\bar{u})-k)-(\bar{u}-k). \\
\end{array}\right.
$$
Since we have $\sgn(x).y\leq |y|$ for any $x,y$ in $\Re^2$, we get out of 
(\ref{edo-proj}) for any $k \in \Re$:
$$
\partial_\tau |\bar{u}-k| \leq \nu\DT\Big(|\Ph(\bar{u})-k|-|\bar{u}-k|\Big).
$$
We observe that $\partial_\tau |\Ph(\bar{u})-k|=0$ along the flow of 
(\ref{edo-proj}) and this implies:
$$
|\bar{u}-k|(x,1)-|\bar{u}-k|(x,0)\leq  \big(1-\exp(-\nu \DT)\big).
\Big(|\Ph(u)-k|-|u-k|\Big)(x,0).
$$
Therefore, we just have to recall the following from (\ref{lim-proj}) to conclude the proof:
$$
\begin{array}{c}
\ds \sum_{n \in \Na_*} \nu \DT \int_\Re \left(\int_0^1 \sgn(\bar{u}-k)
\big(\Ph(\bar{u})-\bar{u}\big)(x,\tau).d\tau\right)\phi(x,n\DT).dx = \\
\ds \sum_{n \in \Na_*} \int_\Re \left(\int_0^1 \partial_\tau |\bar{u}-k|
(x,\tau).d\tau\right)\phi(x,n\DT).dx.
\end{array}
$$
$\Box$

This means in particular that the following inequality holds as a consequence
of (\ref{zozo}) and this last lemma:
\BE \label{zozo-2}
\begin{array}{c}
\ds \Dt \Big [|u-k|+|v-l|\Big] + \Dx \Big |f(u)-f(k) \Big | \leq \\
\ds \sum_{n \in \Na_*} \Big[\big(|\Ph(u)-k|-|u-k|+
|\Ph(v)-l|-|v-l|\big) \big(1-\exp(-\nu \DT)\big) \\
\ds +\mu \DT R(u,v)\big (\sgn(v-l)-\sgn(u-k)\big )\Big]\delta(t-n\DT),\\
\end{array}
\EE
in the sense of Proposition \ref{prop-relax}. One can notice the similarity 
between this last inequality and the one derived in \cite{stw} \S 3; indeed
the numerical approximation considered by these authors matches ours in the
special case where $\nu \rightarrow +\infty$ and a backward Euler solver is
chosen instead of (\ref{edo-relax}). This reduces to a time-splitting scheme
with some kind of ``straight lines" approximation of the quantity 
(\ref{lim-relax}). In contrast, we recall that it has been proved in \cite{bp} 
\S 5, that one has:
$$
\sum_{n \in \Na_*} \nu\DT \Big(\Ph(u)-u\Big)\Big(\sgn(u-k)-\sgn(\Ph(u)-k)\Big).
\delta(t-n\DT) \leq 0,
$$
for any positive continuous test-function in the case where $u \in C^0(\Re_*^+;
L^1(\Re))$.

\subsection{The classical time-splitting scheme}\label{Klassic}

We plan to study the behavior of the system (\ref{systeme3}) in the limit
$\nu \rightarrow +\infty$. Roughly speaking, this reduces to insert a 
projection stage every time $n\DT$, $n \in \Na_*$ {\it before} igniting the 
relaxation mechanism.

\begin{proposition}\label{class-split}
Under the assumptions of Theorem \ref{EU-sing}, the sequence of 
entropy solutions of (\ref{systeme3}) converges strongly in $L^1_{loc}(\Re 
\times \Re_*^+)$ as $\nu \rightarrow +\infty$ towards the  entropy
solution of:
\BE \label{systeme4}
\left\{\begin{array}{c}
\Dt u +\Dx f(u)=-\sum_{n \in \Na_*} \mu \DT R(u,v).\delta(t-n\DT), \\
\Dt v =\sum_{n \in \Na_*} \mu \DT R(u,v).\delta(t-n\DT), \\
\end{array}\right.
\EE
where the right-hand side is still defined by means of (\ref{lim-relax}),
(\ref{edo-relax}), but with the  piecewise constant initial data for $t=n\DT$,
$n \in \Na_*$:
\BE \label{init-relax-split}
\bar{u}(x,\tau=0)=\Ph(u)(x,t-0),\qquad \bar{v}(x,\tau=0)=\Ph(v)(x,t-0).
\EE
Moreover, the following entropy inequality holds true for any nonnegative test
function in ${\cal D}(\Re \times \Re_*^+)$:
\BE \label{zuzu}
\begin{array}{c}
\ds \Dt \Big [|u-k|+|v-l|\Big] + \Dx \Big |f(u)-f(k) \Big | \leq \\
\sum_{n \in \Na_*} \Big\{\big(|\Ph(u)-k|-|u-k|\big) + 
\big(|\Ph(v)-l|-|v-l|\big)+\\ 
\mu \DT.R(u,v)\big(\sgn(v-l)-\sgn(u-k)\big)\Big\}\delta(t-n\DT),
\end{array}
\EE
with $k,l$ in $\Re$ and according to (\ref{lim-relax}), (\ref{edo-relax}),
(\ref{init-relax}).
\end{proposition}

{\bf Proof. }
We split it into several steps for the sake of clarity.

(i) Construction of the solution of (\ref{edo-proj}) by a fixed point argument; for $\tau \in [0,1]$,
we consider any function $w \in C^0(0,1;L^1(\Re))$ together with:
$$
\partial_\tau \bar{u}=\nu \DT (\Ph(w)-\bar{u}), \qquad \bar{u}(.,0) \in L^1(\Re)
$$
whose solution is obvious:
$$
\bar{u}(x,\tau)=\exp(-\nu \DT \tau)\bar{u}(x,0)+\nu \DT\exp(-\nu \DT \tau)
\int_0^\tau \exp(\nu \DT s) \Ph(w)(x,s).ds.
$$
This defines a linear mapping $w \mapsto \bar{u}$ in $C^0(0,1;L^1(\Re))$ for 
any choice of the initial datum. We want to establish a contraction property: 
for all $w_1,w_2$ and any $\tau \in [0,1]$, we observe that
$$
\begin{array}{rcl}
\ds \int_\Re |\bar{u}_1-\bar{u}_2|(\tau).dx &\leq & \nu \DT
\ds \int_0^\tau \exp(\nu \DT (s-\tau)) \int_\Re |\Ph(w_1)-\Ph(w_2)|(x,s).dx.ds \\
& \leq & (1-\exp(-\nu\DT))\|w_1-w_2\|_{C^0(0,1;L^1(\Re))}
\end{array}
$$
By Picard's fixed point theorem, we obtain a unique solution in 
$C^0(0,1;L^1(\Re))$ for any $\nu \in \Re^+$ to:
$$
\partial_\tau \bar{u}=\nu \DT (\Ph(\bar{u})-\bar{u}), \qquad 
\bar{u}(.,0) \in L^1 \cap BV(\Re)
$$
Now, since $\partial_\tau \Ph(\bar{u})=0$ along its flow, this
differential equation therefore admits an explicit solution:
$$
\bar{u}(.,\tau)=\exp(-\nu \DT \tau)\bar{u}(.,0)+
(1-\exp(-\nu \DT \tau))\Ph(\bar{u})(.,0).
$$
In particular, we get:
\BE \label{cons}
\Big(\bar{u}-\Ph(\bar{u})\Big)(.,\tau)=\exp(-\nu \DT \tau)
\Big(\bar{u}-\Ph(\bar{u})\Big)(.,0).
\EE

(ii) Limit as $\nu \rightarrow +\infty$ of the measure  
(\ref{lim-proj}); for any $\phi \in C^0_c(\Re \times \Re^+_*)$, we have:
$$
\begin{array}{c}
\ds \sum_{n \in \Na_*} \nu \DT \int_\Re \left(\int_0^1 
\big(\Ph(\bar{u})-\bar{u}\big)(x,\tau).d \tau \right)\phi(x,n\DT).dx = \\
\ds \sum_{n \in \Na_*} \nu \DT \int_\Re \big(\Ph(u)-u\big)(x,n\DT-0)
\left(\int_0^1 \exp(-\nu \tau) .d \tau \right)\phi(x,n\DT).dx \\
\stackrel{\nu \rightarrow +\infty}{\rightarrow}
\ds \sum_{n \in \Na_*} \int_\Re \big(\Ph(u)-u\big)(x,n\DT-0)
\phi(x,n\DT).dx,
\end{array}
$$
which is a finite sum provided $\DT$ remains strictly positive.

(iii) Strong compactness as $\nu \rightarrow +\infty$: the bound (\ref{equic})
seems to blow up despite the fact that (\ref{Lun}), (\ref{BeVe}) still hold.
Making use of (\ref{cons}) in the last step of the proof of Lemma \ref{compact},
we notice that for any $T \in \Re^+$ and $\xi >0$:
\BE
\begin{array}{c}
\ds \int_{\Re \times [t,t+\xi]} \sum_{n =1}^{1+[T/\DT]} \nu \DT \left\{
\int_0^1 \Big|\bar{u}-\Ph(\bar{u})\Big|+\Big|\bar{v}-\Ph(\bar{v})\Big|(x,\tau)
d\tau \right\}.dx.ds  \leq \\
\xi \frac{h}{\DT} (T + \DT) \Big(TV(u_0)+TV(v_0)\Big)\big(1-\exp(-\nu \DT)\big).
\end{array}
\EE
Therefore, exploiting (\ref{lim-proj}) and Lemma \ref{relaxing}, we can deduce 
a fine estimate for the weak solutions of (\ref{systeme3}), (\ref{lim-proj}), 
(\ref{edo-proj}) which is uniform in both $\nu$ and $\mu$:
\BE \label{equic2}
\begin{array}{c}
\ds \sup_{\xi \not = 0} \int_{\Re \times [0,T]} 
\frac{|u(x,t+\xi)-u(x,t)|}{\xi}+\frac{|v(x,t+\xi)-v(x,t)|}{\xi}.dx.dt
\leq \\
\ds (T+\DT) \Big\{(h/\DT+Lip(f))\Big(TV(u_0)+TV(v_0)\Big) +O(1)\Big\}
\end{array}
\EE
Together with (\ref{Lun}), (\ref{BeVe}), this provides compactness in 
$L^1_{loc}(\Re \times \Re_*^+)$ by means of Helly's theorem.
And in the limit, we get the following system: (recall from (\ref{tpl2}))
$$
\left\{\begin{array}{c}
\Dt u +\Dx f(u)=\sum_{n \in \Na_*}
\Big[(\Ph(u)-u)(.,n\DT-0) - \mu \DT R(u,v)\Big]\delta(t-n\DT), \\
\Dt v =\sum_{n \in \Na_*}
\Big[(\Ph(v)-v)(.,n\DT-0) + \mu \DT R(u,v)\Big]\delta(t-n\DT).
\end{array}\right.
$$
Since we know that $\ae(t)=\be(t-\e)$, we must handle the projection step first
and the initial data for (\ref{edo-relax}) becomes (\ref{init-relax-split}). 
The preceding system then rewrites like (\ref{systeme4}).

(iv) Entropy inequality: Using the same arguments together with Lemma 
\ref{entro-proj}, we derive for any $k \in \Re$:
$$
\begin{array}{c}
\ds \sum_{n \in \Na_*} \nu \DT \int_\Re \left(\int_0^1 \sgn\big(\bar{u}-k\big)
\big(\Ph(\bar{u})-\bar{u}\big)(x,\tau).d \tau \right)\phi(x,n\DT).dx = \\
\ds \sum_{n \in \Na_*} \int_\Re \Big(|\bar{u}-k|(x,1)-|\bar{u}-k|(x,0)\Big)
\phi(x,n\DT).dx \leq \\
\ds \sum_{n \in \Na_*} \int_\Re \Big(|\Ph(u)-k|-|u-k|\Big)(x,n\DT-0)
\big(1-\exp(-\nu \DT)\big)\phi(x,n\DT).dx \\
\stackrel{\nu \rightarrow +\infty}{\rightarrow}
\ds \sum_{n \in \Na_*} \int_\Re \Big(|\Ph(u)-k|-|u-k|\Big)(x,n\DT-0)
\phi(x,n\DT).dx.
\end{array}
$$
This way, we get (\ref{zuzu}) out of (\ref{zozo}) and the contraction property
of Theorem \ref{EU-sing} is also preserved.
$\Box$

One can observe that as soon as the initial data of (\ref{systeme4}),
(\ref{lim-relax}), (\ref{edo-relax}), (\ref{init-relax-split}) is 
$\Ph$--invariant and the CFL condition $\DT \|f'(u)\|_{L^\infty} \leq h$ 
holds, its entropy solution matches the piecewise constant 
approximation generated by a classical time-splitting Godunov scheme,
see {\it e.g.} \cite{areg,stw} and \cite{gtz}. 

\begin{theorem}\label{2}
Under the assumptions (\ref{hypos}), the CFL restriction
$\DT \|f'(u)\|_{L^\infty} = h$ and for $v_0=A(u_0)$, 
$u_0 \in L^1 \cap BV(\Re)$, the entropy solutions of 
(\ref{systeme4}), (\ref{zuzu}), (\ref{lim-relax}), (\ref{edo-relax}), 
(\ref{init-relax-split}) converge strongly in $L^1_{loc}(\Re \times \Re_*^+)$
towards $u$, the entropy solution in the sense of Kru\v{z}kov to:
$$
\Dp{t}{\big(u+A(u)\big)}+\Dp{x}{f(u)}=0, \qquad u(x,0)=u_0.
$$
as $\mu\DT \rightarrow +\infty, h \rightarrow 0$.
\end{theorem}

{\bf Proof. }
As we did before concerning Proposition \ref{class-split}, we are about to split the 
forthcoming proof for the sake of clarity.

(i) Study of the solution of (\ref{edo-relax}); for $\tau \in [0,1]$, we 
rewrite its second equation as:
\BE \label{eqn-mod}
\partial_\tau \bar{v}=\mu \DT (\tilde{A}(\bar{v})-\bar{v}), \qquad
\tilde{A}(\bar{v})=A\big(\bar{u}(x,0)+\bar{v}(x,0)-\bar{v}(x,\tau)\big).
\EE
We follow the classical approach, see {\it e.g.} \cite{bre}, and we introduce
the following Banach space for some $k \in \Re^+$ to be fixed later:
$$
X=\left\{ \varphi \in C^0(0,1;L^1(\Re)) \mbox{ such that } \sup_{t \in [0,1]} 
\Big(\exp(-k.t)\|\varphi(.,t)\|_{L^1(\Re)}\Big) < +\infty \right\}.
$$
For any $w \in X$, we consider the mapping $w \mapsto \bar{v} \in X$,
$$
\bar{v}(x,\tau)=\bar{v}(x,0)\exp(-\mu  \DT \tau)+ 
\int_0^\tau \mu\DT \exp(\mu \DT (s-\tau))\tilde{A}(w)(x,s).ds,
$$
which turns out to be a contraction for $k > \mu\DT( Lip(A)-1)$. By Picard's 
fixed point
theorem, we get a unique solution to (\ref{eqn-mod}) $\bar{v} \in X$ for any 
$\mu \in \Re^+$. Moreover, it satisfies
$$
\partial_\tau (\tilde{A}(\bar{v})-\bar{v})(x,\tau)=-\mu 
\Big(1+A'\big(\bar{u}(x,0)+\bar{v}(x,0)-\bar{v}(x,\tau)\big)\Big)
(\tilde{A}(\bar{v})-\bar{v})(x,\tau),
$$
which leads to:
$$
\big|A(\bar{u})-\bar{v}\big|(x,\tau) \leq \exp(-\mu \DT \tau)
\big|{A}(\bar{u})-\bar{v}\big|(x,0).
$$

(ii) Behavior of (\ref{zuzu}) as $\mu \rightarrow +\infty$; the bound 
(\ref{equic2}) can be simplified by means of the CFL condition:
\BE \label{equic3}
\begin{array}{c}
\ds \sup_{\xi \not = 0} \int_{\Re \times [0,T]} 
\frac{|u(x,t+\xi)-u(x,t)|}{\xi}+\frac{|v(x,t+\xi)-v(x,t)|}{\xi}.dx.dt \leq \\
(T+\DT) \Big\{2Lip(f)\Big(TV(u_0)+TV(v_0)\Big) +O(1)\Big\}.
\end{array}
\EE
This provides strong $L^1_{loc}$ compactness and for any positive
continuous test function $\phi \in C^0_c(\Re \times \Re_*^+)$, we have:
$$
\begin{array}{c}
\ds \sum_{n \in \Na_*} \DT \int_\Re \mu \Big( \int_0^1 
\big(A(\bar{u})-\bar{v}\big) \big(\sgn(\bar{v}-l)-\sgn(\bar{u}-k)\big)
(x,\tau).d\tau \Big)\phi(x,n\DT).dx  \\
\leq 2\ds \sum_{n \in \Na_*} \int_\Re 
\big|A(\Ph(u))-\Ph(v)\big|(x,n\DT)\Big( \int_0^1
\mu \DT .\exp(-\mu \DT \tau).d\tau \Big) \phi(x,n\DT).dx  \\
\leq O\left(\frac{1-\exp(-\mu \DT)}{\mu \DT}\right)\|\phi\|_{C^0}.
\end{array}
$$
We used 
the preceding estimates in conjunction with the
bound in Lemma \ref{relaxing}. Therefore, this term vanishes in the limit 
$\mu\DT \rightarrow +\infty, h \rightarrow 0$. Moreover, by Jensen's inequality 
for convex functions, we also have thanks to (\ref{BeVe}) and the CFL condition,
$$
\begin{array}{c}
\ds \sum_{n \in \Na_*} \int_\Re 
\Big(|\Ph(u)-k|-|u-k|\Big)(x,n\DT-0)\phi(x,n\DT).dx \leq \\
\ds \sum_{n \in \Na_*} \int_\Re 
\big|\Ph(u)-u\big|(x,n\DT-0)\big|\Ph(\phi)-\phi\big|(x,n\DT).dx \leq \\
h \sum_{n \in \Na_*} \DT Lip(f) \big(TV(u_0)+TV(v_0)\big)Lip(\phi)
\stackrel{h \rightarrow 0}{\rightarrow} 0.
\end{array}
$$

(iii) It remains to pass to the limit in (\ref{zuzu}); we perform the following 
splitting for $x \in \Re$ and $t\in [n\DT, (n+1)\DT[$
$$
\begin{array}{rcl}
| v-A(u)|(x,t) &\leq &| v(x,t)-v(x,n\DT)|+|v-\Ph(v)|(x,n\DT)+ \\
&&|A(\Ph(u))-\Ph(v)|(x,n\DT)+ \\
&&|A(\Ph(u))-A(u)|(x,n\DT)+| A(u)(x,t)-A(u)(x,n\DT)|, \\
\end{array}
$$
each term converging to zero in $L^1$ as $\mu\DT \rightarrow +\infty$, $h 
\rightarrow 0$ under the prescribed CFL condition. Therefore, we fix $l=A(k) 
\in \Re$ and we observe that for any positive $\psi \in {\cal D}(\Re \times 
\Re_*^+)$
$$
\int_{\Re \times \Re_*^+} |v-l| \Dt \psi(x,t).dx.dt 
\rightarrow 
\int_{\Re \times \Re_*^+} |A(u)-A(k)| \Dt \psi(x,t).dx.dt,
$$
which leads to the integral form of Kru\v{z}kov's entropy inequality:
$$
\int_{\Re \times \Re_*^+} \Big(|u-k|+|A(u)-A(k)|\Big) \Dt \psi(x,t)
+\Big|f(u)-f(k)\Big| \Dx \psi(x,t).dx.dt \geq 0.$$
$\Box$


\subsection{A modified time-splitting scheme}\label{leVek}
 
It is clear that the estimates of Lemmas \ref{compact} and \ref{relaxing} still 
hold for the slightly modified system (compare with (\ref{systeme2})):
\BE
\label{systeme5}
\left\{\begin{array}{l}
\Dt \ue +\Dx f(\ue)= - \mu R(\ue,\ve)\Dt \ae+ \nu (\Ph(\ue)-\ue)\Dt \be, \\
\Dt \ve = \mu R(\ue,\ve)\Dt \ae + \nu (\Ph(\ve)-\ve)\Dt \be, \\
\ue(.,0)=u_0, \ve(.,0)=v_0.
\end{array}\right.
x \in \Re, t>0
\EE
The Cauchy problem for (\ref{systeme5}) still can be studied within the 
framework of \cite{tv} for any $\e >0$. According to the same ideas, we send
$\e \rightarrow 0$ to get the following result.

\begin{proposition}
Let $(u_0,v_0) \in L^1 \cap BV(\Re)$; under the assumptions (\ref{hypos}), 
the sequence $(\ue,\ve)$ of entropy solutions to (\ref{systeme5}) converges
in $L^1_{loc}(\Re\times \Re_*^+)$ towards the entropy solution in the
sense of Definition \ref{gene-entro-soln} to
\BE
\label{modif-split}
\left\{\begin{array}{c}
\Dt u +\Dx f(u)=\sum_{n \in \Na_*} \DT
\Big[- \mu R(u,v)+\nu (\Ph(u)-u) \Big]\delta(t-n\DT), \\
\Dt v =\sum_{n \in \Na_*} \DT
\Big[\mu R(u,v)+\nu (\Ph(v)-v) \Big]\delta(t-n\DT), \\
\end{array}\right.
\EE
which satisfies for any $k,l \in \Re$ and any nonnegative test function
in ${\cal D}(\Re \times \Re_*^+)$:
\BE \label{zezette}
\begin{array}{c}
\ds \Dt \Big [|u-k|+|v-l|\Big] + \Dx \Big |f(u)-f(k) \Big | \leq \\
\ds \sum_{n \in \Na_*} \DT \Big[\mu R(u,v)\big (\sgn(v-l)-\sgn(u-k)\big ) \\
\ds +\nu (\Ph(u)-u)\sgn(u-k)+\nu (\Ph(v)-v)\sgn(v-l) \Big]\delta(t-n\DT),\\
\end{array}
\EE
in the sense of Propositions \ref{prop-proj} and \ref{prop-relax}. Given 
any other set of initial data $(\tilde{u}_0,\tilde{v}_0)$ in $L^1 \cap BV(\Re)$,
this last inequality implies the contraction property for any $t \in \Re^+$:
$$
\|u(.,t)-\tilde{u}(.,t)\|_{L^1(\Re)}+\|v(.,t)-\tilde{v}(.,t)\|_{L^1(\Re)} \leq 
\|u_0-\tilde{u}_0\|_{L^1(\Re)}+\|v_0-\tilde{v}_0\|_{L^1(\Re)}.
$$
\end{proposition}

{\bf Proof. }
Similar estimates than the ones shown in the proof of Lemma \ref{compact} 
ensure the strong $L^1_{loc}$ compactness. The inequality (\ref{zezette})
comes from (\ref{sta}) together with Propositions \ref{prop-proj} and 
\ref{prop-relax}. The contraction property in $L^1(\Re)$ is proved the
same way than for Theorem \ref{EU-sing}.
$\Box$

\begin{remark}
It is important to notice that the entropy solutions of 
(\ref{modif-split}) do {\bf not} coincide with the ones of (\ref{systeme3})
since the relaxation and the projection steps are not applied in
the same order. This means in particular that one has to be careful
when dealing with the right-hand sides as different processes
are concentrated inside the Dirac masses. Hence there
is no uniqueness problem for both systems (\ref{systeme3}), 
(\ref{modif-split}) relying on Propositions \ref{prop-proj} and 
\ref{prop-relax} together with the inequalities (\ref{zozo}) and 
(\ref{zezette}).
\end{remark}

We plan now to carry out the same program as we did for the classical 
time-splitting scheme in the preceding section. By letting $\nu \rightarrow
+ \infty$, we derive an equation whose entropy solution coincide
with the piecewise constant approximation generated by a numerical scheme 
designed in the context of reactive Euler equations for which the reaction
process operates {\it before} the projection stage, \cite{lev}.

\begin{proposition}
Assume (\ref{hypos}) and $(u_0,v_0) \in L^1 \cap BV(\Re)$, the sequence of 
entropy solutions to (\ref{modif-split}) converges strongly in
$L^1_{loc}(\Re \times \Re_*^+)$ as $\nu \rightarrow +\infty$ towards the
entropy solution of:
\BE\label{modif-split-2}
\left.\begin{array}{c}
\Dt u +\Dx f(u)=\sum_{n \in \Na_*}
\Big[- \mu\DT R(u,v)+(\Ph(u)-u)(.,n\DT+0) \Big]\delta(t-n\DT), \\
\Dt v =\sum_{n \in \Na_*} 
\Big[\mu\DT R(u,v)+(\Ph(v)-v)(.,n\DT+0) \Big]\delta(t-n\DT), \\
\end{array}\right.
\EE
where the right-hand sides are defined by means of (\ref{lim-relax}), 
(\ref{edo-relax}), (\ref{init-relax}).
Moreover, the following entropy inequality holds true for any nonnegative test
function in ${\cal D}(\Re \times \Re_*^+)$ with $k,l$ in $\Re$:
\BE \label{zuzu-bis}
\begin{array}{c}
\ds \Dt \Big [|u-k|+|v-l|\Big] + \Dx \Big |f(u)-f(k) \Big | \leq \\
\sum_{n \in \Na_*} \Big\{\mu\DT.R(u,v)\big(\sgn(v-l)-\sgn(u-k)\big)+ \\
\big(|\Ph(u)-k|-|u-k|\big) + 
\big(|\Ph(v)-l|-|v-l|\big) \Big\}\delta(t-n\DT),
\end{array}
\EE
\end{proposition}
The final statement is concerned with the relaxation limit in 
(\ref{modif-split-2}).

\begin{theorem}\label{4}
Under the assumptions (\ref{hypos}), the CFL restriction
$\DT \|f'(u)\|_{L^\infty} = h$ and for $v_0=A(u_0)$,
$u_0 \in L^1 \cap BV(\Re)$, the entropy solutions of
(\ref{modif-split-2}), (\ref{zuzu-bis}), (\ref{lim-relax}), (\ref{edo-relax}), 
(\ref{init-relax}) converge strongly in $L^1_{loc}(\Re \times \Re_*^+)$
towards $u$, the entropy solution in the sense of Kru\v{z}kov to:
$$
\Dp{t}{\big(u+A(u)\big)}+\Dp{x}{f(u)}=0, \qquad u(x,0)=u_0.
$$
as $\mu\DT \rightarrow +\infty, h \rightarrow 0$.
\end{theorem}

The proofs of these last two statements are completely similar to the ones of
Proposition \ref{class-split} and Theorem \ref{2} hence we skip them.

\section{Conclusion} 

In this paper, we considered some ``measure source terms" for a quasilinear relaxation
system whose weak solutions coincide with piecewise constant approximations generated
by commonly-used time-splitting numerical schemes. This allows to establish
convergence results directly relying on properties of the underlying original
system and bypassing some heavy computations required to show the 
stability of the numerical processes (see {\it e.g.} \S 2.5 in \cite{moc}). 

\paragraph{Acknowledgments: } The author would like to thank Professor
Wen-An Yong for several valuable discussions and suggestions. This work has been
partially supported by the European Union TMR project  HCL \#ERBFMRXCT960033.

\nonumsection{References}

\end{document}